# EXPLORING COLLABORATION: THE EFFECT OF GENDER ON MATHEMATICS LEARNING PREFERENCES

Sang Hyun Kim and Tanya Evans

University of Auckland

*This study examines the influence of gender on students' collaborative preferences for learning mathematics (CPLM) over time in an undergraduate mathematics context. Data collected at three points during the semester were analyzed using a two-way mixed ANOVA. Results showed no significant interaction between gender and time, nor a main effect of gender, indicating stable CPLM scores and comparable preferences between male and female students. These findings challenge prior research on gendered differences in collaboration and suggest that further exploration of contextual factors is needed to deepen understanding of CPLM in undergraduate mathematics education.*

## INTRODUCTION

Amid the growing interest in improving student engagement within undergraduate mathematics education, understanding the diverse preferences that students bring to their learning has become increasingly important. These preferences, particularly in collaborative settings, are shaped by a complex interplay of cognitive, affective, and social factors, and it seems likely that gender may emerge as a salient factor that affects them. The persistent gender differences reported in the literature may suggest that gender can shape students' preferences to collaborate with their peers in mathematics learning contexts (e.g., Else-Quest et al., 2010; Seegers & Boekaerts, 1996). Exploring these dynamics can provide valuable insights into whether gender differences are present and how they may influence student preferences and behaviors within an undergraduate mathematics context. This study seeks to contribute to this area of research by examining the role of gender in shaping collaborative preferences for learning mathematics.

The importance of understanding such collaborative preferences of students is that it may provide us with a greater understanding of how the dispositions held by an individual may be linked to other related social-cognitive and affective factors. In certain contexts, stronger preferences to engage with mathematics in a collaborative manner are likely to drive behaviors within and beyond the classroom that align with them, which could impact learning outcomes, including achievement. A preference for collaboration may foster active participation and idea sharing, promoting social learning processes that enhance conceptual understanding. Alternatively, preferences for independent and individualized work might encourage one to work at a more appropriate pace and develop greater autonomy while problem-solving. Furthermore, these preferences may manifest in contexts beyond academic settings, influencing how students engage with others and apply mathematical reasoning to real-world problems. In this way, preferences act as a lens through which we can better understand the interplay between students' social tendencies, affective experiences, and their practical use of mathematical knowledge.

## Preferences in mathematics

Student preferences in mathematics education reflect individual inclinations toward specific modes of engagement, such as collaborative or independent work. Preferences are an example of a motivational structure that serves as "templates for whether to put forth effort towards mathematical activity and the extent to which efforts are seen as efficacious" (Goldin et al., 2016, p. 18). These preferences are closely tied to students' values and can help shape the tasks they engage with and the pedagogical approach in which they do so. Beliefs are typically grounded in students' perceptions of truth or correctness (Beswick, 2007). On the other hand, preferences are shaped by personal experiences and interests, influencing how students engage with mathematical tasks in various contexts. While preferences can guide short-term behaviors, such as the choice to work with others or tackle a problem independently, they can also evolve into more stable patterns or habits over time, influencing long-term academic engagement (Goldin et al., 2016). These stable preferences can impact not only students' mathematical achievement but also their persistence and development of mathematical identities. Hence, understanding the role of preferences may be important for informing teaching practices that align with students' inclinations and fostering more effective learning experiences.

## Collaborative preferences

A desire to work more collaboratively or individually is a natural point of variation worth exploring, with limited studies exploring this at the undergraduate level or within the domain of mathematics. The literature on student experiences during collaboration underscores a strong need to understand ways in which we can better capture the preferences and attitudes of students towards collaboration within mathematics. Existing studies in the broader educational literature have directed attention to some factors that may be correlated with certain preferences. Gajderowicz et al. (2023) found that students preferred working in groups over working individually, although these preferences for group work depended on factors such as subject domain and ability. Furthermore, the authors concluded that understanding student preferences can improve student satisfaction in learning mathematics, noting that "students can enjoy learning math if comfortable study modes are used" (Gajderowicz et al., 2023, pp. 8-9), a conclusion supported by earlier STEM education research (Okebukola, 1986). In fact, Okebukola (1986) found that alignment between students' preferred modes of learning and the mode in which they were taught resulted in the best achievement outcomes. Similarly, aligning teaching practices with student preferences for engagement can result in improvements in conceptual learning and perceptions of autonomy, as demonstrated by Jang et al. (2016). However, some scholars warn against simply following the desires of the students, as while it could create an illusion of learning, this may not always be reflective of an improvement in student achievement (Deslauriers et al., 2019).

## Gender differences

Research on gender differences in mathematics achievement has shown mixed findings, with boys historically outperforming girls in certain contexts, particularly at higher levels of mathematics and in specific regions (Else-Quest et al., 2010; Guiso et al., 2008; Seegers & Boekaerts, 1996); however, this gap has diminished or disappeared in many countries over recent decades (Hyde & Mertz, 2009). Some differences between males and females in affective and motivational factors and their effects on mathematics achievement, including self-concept, attribution of success, and perceptions of difficulty, have also been reported (Ethington, 1992; Mejía-Rodríguez et al., 2021; Wolleat et al., 1980). Research suggests that girls often prefer collaborative environments emphasizing collective reasoning and valuing relationships, while boys may be more inclined to work individually (Gajderowicz et al., 2023; Kanevsky et al., 2022; Owens, 1985). Borgonovi et al. (2023) further contribute to this understanding by emphasizing that these gendered preferences are not solely reflective of individual inclinations but rather are influenced by a range of cognitive, emotional, and social factors. Using 2015 PISA data, the authors reported that girls scored higher than boys in collaborative problem-solving skills across all countries, with this difference being more pronounced in favour of girls in countries with greater gender equality (Borgonovi et al., 2023). These findings highlight how societal contexts and cultural norms shape not only achievement gaps but also students' preferences and approaches to learning. This insight underscores the need to consider gender as a dynamic and influential factor in shaping students' collaborative preferences for learning mathematics within the context of undergraduate education while recognizing that collaboration and student engagement may evolve over time.

## RESEARCH AIM

The aim of this study is to explore the influence of gender on students' collaborative preferences for learning mathematics (CPLM) over time. This research addresses two primary questions: (1) *what are the collaborative preferences for learning mathematics among undergraduate students across different genders*, and (2) *how does gender (male vs. female) influence students' collaborative preferences for learning mathematics (CPLM) across different time points?*

## METHODOLOGY

### Context

This study draws on data gathered during the second semester of 2024 from a second-year undergraduate service mathematics course at a large university in New Zealand. This course introduces students to key mathematical concepts required for various academic programs across the university, covering three main topics: calculus II, linear algebra II, and differential equations. Delivered over a 12-week semester, the course included weekly one-hour problem-solving sessions (tutorials) in which students were given a set of mathematics problems and were encouraged to work with peers, although this was not mandated. A tutor (either the lecturer or an experienced graduate student) would be available to assist students with problem-solving by providing guidance, clarifying concepts, and answering questions as needed.

**Data**

The self-report data used in this study was part of a larger project, with data collected at three points in the semester (the first week, after the midsemester break, and the final week). Here, we focus on students' collaborative preferences for learning mathematics (CPLM) and gender. Of the 294 students enrolled in the course, 201 completed the surveys at all three time points. Students who selected 'declined to answer' about their gender (n = 3) were not included in the analysis, resulting in a sample size of 198 students for the analysis.

Students' collaborative preferences for learning mathematics were assessed using a 5-item scale developed and validated by Kim and Evans (2025). Students responded to five slider-based, close-ended items, rating their preferences on a scale from 0 (indicating a preference for individual learning) to 100 (indicating a preference for collaborative learning). These ratings captured students' collaborative preferences across various scenarios. Two example items from the scale are *"What is the most effective way for you to learn mathematics?" and "In what social setting do you prefer to be exposed to novel concepts?"* To assess students' collaborative preferences for learning mathematics (CPLM), we calculated a composite score for each student at each time point by averaging their responses across the five items.

**Analysis**

Within this study, we used descriptive statistics to examine students' CPLM based on gender. We also conducted a two-way mixed ANOVA to investigate the effect of gender as a between-subject factor on CPLM scores over time measured at three time points. Gender was the between-subjects factor (male and female) and time was the within-subjects factor (3 levels), with an interaction term to assess whether the change in CPLM scores over time differed by gender.

**RESULTS**

Table 1 shows the descriptive statistics for the CPLM scores for males and females at the three time points. As shown in Figure 1, the CPLM scores for males were slightly higher than those for females at each time point, suggesting a higher preference for collaboration. The differences between the means of males and females also did not differ significantly across the semester, with the differences between male and female students' CPLM scores at the three time points being 4.44, 3.16, and 5.23, respectively.

| CPLM | Gender | Mean | Std. Deviation | N |
|------|--------|------|----------------|---|
| Time 1 | Male | 50.40 | 22.34 | 115 |
| | Female | 45.96 | 22.28 | 83 |
| | Total | 48.54 | 22.37 | 198 |
| Time 2 | Male | 49.19 | 23.13 | 115 |
| | Female | 46.03 | 21.73 | 83 |
| | Total | 47.86 | 22.55 | 198 |
| Time 3 | Male | 52.19 | 23.73 | 115 |
| | Female | 46.96 | 21.94 | 83 |
| | Total | 50.00 | 23.08 | 198 |

Table 1: Descriptive statistics for CPLM data by gender across three time points during the semester.

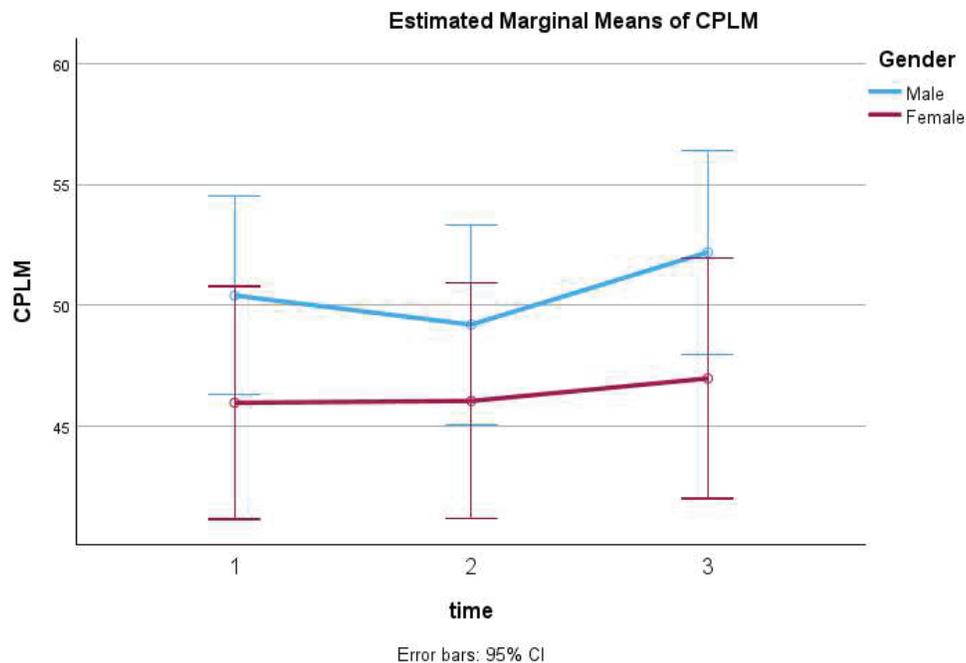

Figure 1: Change in CPLM scores over time for males and females.

Mauchly's test of sphericity revealed a violation of the sphericity assumption for the two-way interaction, $\chi^2(2) = 12.89$, $p = .002$. Greenhouse-Geisser corrections were applied to adjust the degrees of freedom for the $F$-tests going forward with the analysis of within-subject effects in order to account for this violation ($\varepsilon = .904$). No statistically significant interaction between gender and time on CPLM was found, $F(1.880, 368.429) = .260$, $p = .758$, partial $\eta^2 = .001$. The main effect of time did not show a statistically significant difference in mean CPLM at the different time points, $F(1.880, 368.429) = .976$, $p = .373$, partial $\eta^2 = .005$.

Similarly, the main effect of gender did not show a statistically significant difference in mean CPLM between males and females, $F(1, 196) = 2.340$, $p = .128$, partial $\eta^2 = .012$. While these differences were not statistically significant ($p = .128$), the partial eta-squared value ($\eta^2 = .012$) suggests a small effect size, indicating the potential for a gender-related difference. This may warrant further investigation, ideally with a larger sample or over an extended time frame.

## DISCUSSION

This study investigated collaborative preferences for learning mathematics (CPLM) within a sample of undergraduate students and whether gender influenced these preferences over time. The mean scores for males were higher than for females at all time points, with the difference between them largest by the end of the semester (5.23). While this suggests males exhibited a slightly higher preference for collaborative learning than females, contrary to prior research findings, these differences were not statistically significant. Notably, the set of 95% confidence intervals of mean values across genders and time points, ranging from 41.1 to 56.4, suggest that most students maintained a balanced preference for both collaborative and individual learning.

Interestingly, the existing literature typically associates males with a preference for individual work (e.g., Gajderowicz et al., 2023; Owens, 1985). Such findings tend to emphasize the social and emotional aspects of collaborative learning, which are traditionally seen as aligning more closely with gender norms favoring girls. However, our results challenge the existing literature, as gender differences in CPLM were not only statistically non-significant but also reversed, with boys preferring more collaboration than girls. This trend was consistent across all three time points during the semester. This could reflect contextual factors specific to this study, such as the structure of the undergraduate tutorials, where collaboration was encouraged but not essential, and where marks were awarded for participation regardless of the work done individually or collaboratively, potentially mitigating any effect on student preferences. Another possible explanation, from a psychological perspective, is the potential benefits that male students may derive from this type of tutorial environment. Raabe and Block (2024) suggest that males may derive social validation and confidence from their performance in mathematics, which is reinforced through peer interactions and friendship networks, whereas females do not feel the same social pressure.

Lastly, no significant interaction effect between gender and time was observed in this study, nor was the main effect of gender found on collaborative preferences for learning mathematics (CPLM). These findings suggest that, within this sample, male and female students exhibited more comparable preferences for collaboration than expected. Moreover, the absence of the main effect of time supports the possibility that CPLM preferences are relatively stable, indicating it may possess a trait-like quality. This stability could stem from broader and more diverse learning experiences undergraduate students have encountered, both in mathematics education and education more generally, compared to younger students. Regardless, the underlying mechanisms are yet to be understood, warranting further research to examine how transient these preferences are and the influence of social factors around them.

This study has several limitations that future research could address. The small, context-specific sample may limit the generalizability of these findings, emphasizing the need for larger, more diverse samples to better assess CPLM stability and detect subtle gender differences with greater statistical power. In particular, the absence of significant gender differences in this study could be due to the sample size, suggesting that a larger sample may be needed in future studies to reduce the likelihood of a Type II error. Additionally, the study did not consider factors such as cultural background, prior collaborative experiences, or academic achievement, all of which could influence CPLM. Future research could address these gaps by incorporating diverse samples and examining additional factors. Furthermore, integrating observational or interview data could provide deeper insights into how students work and how their preferences develop or change in response to specific events, tasks, or interventions.

In conclusion, this research sought to contribute to our understanding of how collaborative preferences and gender intersect in undergraduate mathematics. Future research can build upon this work by continuing to examine the underlying mechanisms behind these preferences, further uncovering the factors that shape collaborative preferences and their relationship with gender in educational contexts.